 \documentclass[12pt,reqno]{amsart}
\usepackage{amsfonts}
\usepackage{amssymb}
\usepackage{amsthm}
\usepackage{upref}
\usepackage{enumerate}

\usepackage{amscd}

\makeatletter 

\def\LaTeX{\leavevmode L\raise.42ex
    \hbox{\kern-.3em\size{\sf@size}{0pt}\selectfont A}\kern-.15em\TeX}
 \makeatother
 
\DeclareMathOperator{\clos}{clos}

\numberwithin{equation}{section}

\newtheorem{lemma}{Lemma}[section]
\newtheorem{theorem}[lemma]{Theorem} 
\newtheorem{corollary}[lemma]{Corollary}
\newtheorem{proposition}[lemma]{Proposition}

\theoremstyle{definition}

\newtheorem{remark}[lemma]{Remark}

  \newcommand{\e}{\eqref}

\newcommand{\q}{\quad}

\newcommand{\ov}{\overline}
\newcommand{\wt}{\widetilde}

\newcommand{\ti}{\tilde}

\newcommand{\la}{\langle}
\newcommand{\ra}{\rangle}

\renewcommand\Im{\operatorname{Im}}

\newenvironment{pf}{\begin{proof}}{\end{proof}}

\def\qqq{\mathrel{\subset\mkern-15mu\lower.38ex\hbox{${\scriptscriptstyle\rightarrow}$}}}

\let\cal\mathcal

\let\Bbb\mathbb
 
 \UseRawInputEncoding 
 
\begin{document}
\title 
{Semiclassical analysis  in the limit circle case}
\author{ D. R. Yafaev  }
\address{   Univ  Rennes, CNRS, IRMAR-UMR 6625, F-35000
    Rennes, France, SPGU, Univ. Nab. 7/9, Saint Petersburg, 199034 Russia, and     NTU Sirius, Olympiysky av. 1, Sochi, 354340 Russia}
\email{yafaev@univ-rennes1.fr}


\thanks {Supported by  RFBR grant No.   20-01-00451 A}

\begin{abstract}
We consider    second order differential  equations   with real coefficients   that are  in the limit circle case at infinity.   Using the semiclassical Ansatz,  we construct solutions (the Jost solutions) of such equations with  a prescribed asymptotic behavior  
for $x\to\infty$.  It turns out that   in the limit circle case,  this Ansatz can be chosen common for all values
of the spectral parameter $z$.  This leads to asymptotic formulas for all solutions   of considered  differential equations, both homogeneous and non-homogeneous. We also efficiently describe all self-adjoint realizations of the corresponding differential operators in terms of boundary conditions at infinity and   
  find a representation  for their resolvents.  
 \end{abstract}

   \dedicatory{To the memory of Serezha Naboko}

 \maketitle

\section{Introduction}

 \subsection{Setting the problem}

We consider a second order differential equation  
\begin{equation}
- (p(x) u' (x))'  - q(x) u(x)=z u(x), \q p(x)>0, \; q(x)=\bar{q} (x),\q x \in {\Bbb R}_{+}, \; z \in {\Bbb C},
  \label{eq:Jy}\end{equation}
  under some minimal regularity assumptions on coefficients $p(x)$ and $q(x)$ guaranteeing that  its solutions, as well as their derivatives, are continuous on $[0,\infty)$.   We are interested in the limit circle  (LC) case at infinity
   where all  solutions are in $L^2 ({\Bbb R}_{+})$ for all $z \in {\Bbb C}$.
   Equation  \e{eq:Jy}  is known as the Schr\"odinger  equation if $p(x)=1$, but we keep the same term in the general case. We write $-q(x)$ because under our assumptions $q(x)\to+ \infty$ as $x\to\infty$.
    
    Our analysis relies on a construction of solutions $f_{z} (x)$ of equation \e{eq:Jy}  distinguished by their asymptotics
 \begin{equation}
 f_{z}  (x)= (p(x) q(x))^{-1/4}  e^{i\Xi (x)} \big( 1+ o(1)\big) , \q z\in {\Bbb C},\q  x\to\infty,
  \label{eq:LG}\end{equation}
  where
  \[
\Xi (x)= \int_{x_{0}}^x \big(q(y)/p(y)\big)^{1/2} dy
\]
($x_{0}$  is an arbitrary fixed number). Solutions    $f_{z} (x)$ are known as the Jost solutions of equation \e{eq:Jy}.  Relation \e{eq:LG} shows that the leading terms of asymptotics of the Jost solutions do not depend on $z\in {\Bbb C}$. This fact is specific for the LC case. The  functions  $f_{z} (x)$ and  $\bar{f}_{\bar{z} }(x)$   satisfy the same equation \e{eq:Jy} and  are linearly independent so that an arbitrary solution of equation \e{eq:Jy}  is their linear combination.
          We require that
     \[
     \int_{x_{0}}^\infty (p(x) q(x))^{-1/2}  dx <\infty
     \]
     whence $ f_{z} \in L^2 ({\Bbb R}_{+})$ for all $z\in{\Bbb C}$. Thus,  according to \e{eq:LG}, we are in the LC case.
     

    \subsection{Limit point versus limit circle}

     The Weyl limit point/circle theory  states that   differential   equation \e{eq:Jy} always has a non-trivial solution in $L^2 ({\Bbb R}_{+})$ for $\Im z \neq 0$. This solution is either unique (up to a constant factor) or all solutions of \e{eq:Jy} belong to $L^2 ({\Bbb R}_{+})$.
The first instance is known as the limit point (LP) case and the second one --  as the limit circle (LC) case. Consistent presentations of the Weyl theory can be found, for example,  in the books \cite{CoLe}, Chapter~IX,  \cite{DuSc},  Chapter~XIII and \cite{RS},  Chapter~X.1.

Our goal is to study  self-adjoint operators  in 
 the space $L^2 ({\Bbb R}_{+})$ associated with a differential operator
    \begin{equation}
(  {\sf H} u)(x)=- \big(p(x) u' (x)\big)'  - q(x) u(x)   .
  \label{eq:ZP+}\end{equation}
  The operator ${\sf H}$ defined on a domain $C_{0} ^\infty ({\Bbb R}_{+}) $ is symmetric in $L^2 ({\Bbb R}_{+})$, but to make it self-adjoint, one has to add boundary  conditions at $x=0$ and, eventually, for $x\to\infty$.  The boundary condition at the point $x=0$  looks as 
     \begin{equation}
u'(0)= \alpha  u(0) \q\mbox{where}\q \alpha =\bar{\alpha}.
  \label{eq:BCz}\end{equation}
  The value $\alpha=\infty$ is not excluded. In this  case \e{eq:BCz} should be understood as the equality $u(0)=0$. We always require   condition \e{eq:BCz}, fix $\alpha$ and   do not  keep track of $\alpha$ in notation. 
  
  Let us define    a symmetric operator $H_{00}$   by an equality $H_{00}u= {\sf H} u $
 on a domain     $  {\cal D} (H_{00}) $ consisting of smooth functions $u(x)$  satisfying    boundary condition \e{eq:BCz} and such that $u(x)=0$  for sufficiently large $x$. 
 The operator $H_{00}$ is   essentially self-adjoint if and only if the  LP case occurs (see, e.g., Theorem~X.7 in the book \cite{RS}). 
  In the LC case,  $H_{00}$ has  a one parameter family $H_{\omega}$ (here $\omega$ is a point on the unit circle ${\Bbb T} \subset {\Bbb C} $) of self-adjoint extensions  distinguished by some   conditions for $x\to\infty$.  Their description can be performed in different terms. Our analysis  relies on asymptotic formula \e{eq:LG}.  
  Alternative possibilities are briefly discussed in Sect.~4.5.

  \subsection{Plan of the paper}
     
     The existence of the Jost solutions is proven in Sect.~2. Another important element of our approach is
     a construction of an operator ${\cal R} (z)$ (we call it the quasiresolvent) playing the role of the resolvent for the maximal operator $H_{\rm max}:= H_{00}^*$ (the adjoint of $H_{00}$). We emphasize that ${\cal R} (z)$ is an analytic function of $z\in{\Bbb C}$.  The operator ${\cal R} (z)$ is constructed in Sect.~3. This section is close to 
     Sect.~2 of paper \cite{Schr-LC}. 
     
 
 In Sect.~4 we show (see Theorem~\ref{LCR}) that all functions   $  u\in {\cal D} (H_{\rm max}) $ have asymptotic behavior
  \begin{equation}
u (x) = (p(x) q(x))^{-1/4} \big(s_{+}   e^{i \Xi (x)} +s_{-}  e^{-i  \Xi (x)}+ o(1)\big)      , \q x\to\infty,
\label{eq:LC2}\end{equation}
with some coefficients $s_{\pm}= s_{\pm} (u)\in {\Bbb C}$. We distinguish a set $  {\cal D} (H_\omega) \subset {\cal D} (H_{\rm max}) $   by  the condition
  \begin{equation}
   s_{+} (u) =\omega s_{-} (u) \q \mbox{where} \q |\omega|=1 .
      \label{eq:gen2}\end{equation}
      Asymptotic coefficients $s_{\pm} (u)$ in  formula \e{eq:LC2} play the role of  boundary values $u(0)$ and 
  $u'(0)$ for functions $u$ in the local Sobolev space ${\Bbb H}^2_{\rm loc}$, and equality \e{eq:gen2} plays the role of  boundary condition \e{eq:BCz}. 
  
        Theorem~\ref{ABX}  shows that the restriction $H_{\omega}$ of the operator  $H_{\rm max}$ on $  {\cal D} (H_\omega)$ is self-adjoint, and all self-adjoint extensions of the operator  $H_{\rm min}$ coincide with one of the operators $H_{\omega}$.  Our proofs of these results are independent of the von Neumann formulas.
      Finally, we construct    the resolvents of the operators $H_{\omega}$ in Theorem~\ref{ABX1}.

   The construction of this paper is similar to the approach  used in Sect.~3 of \cite{Jacobi-LC} in the case of Jacobi operators.  Conditions on coefficients  look completely differently for  Jacobi and Schr\"odinger operators, but,  in both cases, we are in the LC case and the leading terms of asymptotics of the corresponding Jost solutions  do not depend on the spectral parameter. This implies that spectral properties of these  two classes of  operators are similar.



  \section{The semiclassical Ansatz}
  
  We here construct  solutions of the Schr\"odinger equation \e{eq:Jy} with asymptotics \e{eq:LG}  for $x\to\infty$.

 \subsection{Regular solutions}


    To avoid inessential technical complications, we
    always suppose that $ p \in  C^1 ({\Bbb R}_{+})$, $ q \in  C ({\Bbb R}_{+})$ and the functions $p(x)$,   $q(x)$ have finite limits as $x\to 0$. We assume that $p(x)> 0$  for $x\geq 0$. The solutions of equation \e{eq:Jy} exist, belong to $C^2 ({\Bbb R}_{+})$ and they have limits $u(+0)= : u(0)$, $u'(+0)= : u'(0)$. A  solution $u(x)$ is distinguished uniquely by  boundary conditions $u(0)=u_{0}$, $u'(0)=u_{1}$.

 Recall that for arbitrary solutions $u$ and $v$ of equation \e{eq:Jy}  their Wronskian
  \[
\{ u, v  \} : = p(x)   (u' (x)  v (x)- u (x)  v' (x))
\]
does not depend on $x\in{\Bbb R}_{+}$. Clearly, the Wronskian $\{ u, v  \} =0$ if and only if the solutions $u$ and $v$ are proportional.


 We introduce a couple of regular  solutions of   equation \e{eq:Jy} by boundary conditions at the point $x=0$:
 \begin{equation}
  \begin{cases}
\varphi_{z}(0)= 1,&\q  \varphi'_{z}(0)=\alpha, 
\\
\theta_{z}(0)= 0,&\q \theta'_{z}(0)=- p(0)^{-1}, 
 \end{cases} \q\mbox{if}\q \alpha\in{\Bbb R}
  \label{eq:Pz}\end{equation} 
 and
   \begin{equation}
  \begin{cases}
\varphi_{z}(0)= 0,&\q  \varphi'_{z}(0)=1, 
\\
\theta_{z}(0)= p(0)^{-1},&\q \theta'_{z}(0)= 0, 
 \end{cases}   \q\mbox{if}\q \alpha=\infty.
\label{eq:Pzi}\end{equation}
 Obviously,   $\varphi_{z} (x)$ (but not $\theta_{z} (x)$) satisfy boundary condition \e{eq:BCz} and   the Wronskian
   $\{\varphi_{z}, \theta_{z}\}=1$.

In the LC case all solutions of equation \e{eq:Jy} are in  $ L^2 ({\Bbb R}_{+})$. In particular,
 \begin{equation}
\varphi_{z}\in L^2 ({\Bbb R}_{+}),\q \theta_{z}\in L^2 ({\Bbb R}_{+})\q \mbox{for all}\q z\in {\Bbb C}.
  \label{eq:PQz}\end{equation} 
   
  \subsection{Jost solutions}
  
  The Jost solutions $f_{z} (x)$ of  the differential equation \e{eq:Jy}
  are distinguished            by their asymptotics  \e{eq:LG} for $x\to\infty$.     They will be constructed in Theorem~\ref{JOST}.   Let us set
  \begin{equation}
a(x)= (p(x) q(x))^{-1/4},\q \xi (x)= \sqrt{\frac{q(x)}{p(x)}}\q \mbox{and}\q \Xi (x)=\int_{x_{0}}^x \xi(y) dy
  \label{eq:Jo1}\end{equation}
  so that
  \begin{equation}
  p(x) \xi (x) a^2 (x)=1. 
  \label{eq:Jo1A}\end{equation}
  Notation \e{eq:Jo1} will be used throughout the whole paper.
  
  \begin{theorem}\label{JOST}
    Suppose that
      \begin{equation}
a\in L^2 (x_{0}, \infty) \q \mbox{and}\q a(pa')'\in L^1 (x_{0}, \infty)
  \label{eq:JJ}\end{equation} 
  for  some $x_{0}> 0$.
  Then for all $z \in {\Bbb C}$, equation \e{eq:Jy} has a solution  $f_{z} (x)$ with asymptotics
  \begin{equation}
f_{z} (x)= a(x) e^{i\Xi (x)}(1+ o(1))
\label{eq:Jo2}\end{equation}
as $x\to\infty$.  If, additionally,
      \begin{equation}
\lim_{x\to\infty} p(x)  a' (x) a(x)=0,
  \label{eq:JJ+}\end{equation} 
  then
    \begin{equation}
f_{z}' (x)= i \xi(x)a(x) e^{i\Xi (x)}(1+ o(1)), \q x\to\infty.
\label{eq:Jo2+}\end{equation}
   \end{theorem} 
   
     \begin{remark}\label{JOSTr}
     In the leading particular case  $ p (x)=1$ conditions  \e{eq:JJ}  mean that
   \[
   \int_{x_{0}}^\infty q(x)^{-1/2}dx <\infty  \q \mbox{and}\q    \int_{x_{0}}^\infty q(x)^{-1/4}
   | (q(x)^{-1/4})''  | dx <\infty .
\]
  Condition  \e{eq:JJ+}  reduces to $q(x)^{-3/2}
    q'(x)                                                                                                                                                                                              = o(1)$ for $x\to\infty$.
        \end{remark}

 In a detailed notation, formula   \e{eq:Jo2} coincides with \e{eq:LG}.
We call   $f_{z}(x)$ the Jost solution. We emphasize that the leading term of asymptotics of $f_{z}(x)$ does not depend on $z$. This is specific for the LC case.  Otherwise a proof of Theorem~\ref{JOST}  is relatively standard (cf. \cite{Olver}, Chapter~6).
It relies on the fact that the Liouville-Green {\it Ansatz} 
 \begin{equation}
A(x)= a(x) e^{i \Xi(x)}
  \label{eq:Jo3}\end{equation}
  satisfies equation \e{eq:Jy} with a sufficiently good accuracy.
  
  We start a proof of Theorem~\ref{JOST} with  a multiplicative change of variables
  \begin{equation}
f_{z} (x)= A(x)  \psi_{z} (x).
\label{eq:Jo4}\end{equation}
The following result will be obtained by a direct computation.
  
  \begin{lemma}\label{JOST1}
  Set 
   \begin{equation}
\rho_{z}(x)= a(x) \big( p(x) a'(x))'+ z a^2 (x).
  \label{eq:Jr}\end{equation} 
    Then the  equation
    \begin{equation}
 (p(x) f_{z}' (x))'  + q(x) f_{z}(x) + z f_{z}(x)  =0
  \label{eq:JY}\end{equation} 
    is equivalent to an equation
  \begin{equation}
\psi_{z}'' (x) +\big(2i\xi(x) -\frac{\xi'(x)}{\xi(x)}  \big) \psi_{z}'(x)+ \rho_{z}(x)\xi(x) \psi_{z} (x)=0
\label{eq:Jo5}\end{equation}
for the function $\psi_{z} (x)$ defined by formula \e{eq:Jo4}.
    \end{lemma} 
    
    \begin{pf}
    Differentiating \e{eq:Jo4} twice, we find that
    \[
    (pf_{z}')'= pA\psi_{z}'' + (2p A'+ p' A) \psi_{z}'+ (pA')'\psi_{z}.
    \]
    Substituting this expression into \e{eq:JY} and dividing by $pA$, we rewrite \e{eq:JY} as 
     \begin{equation}
\psi_{z}'' +\big(2\frac{A'}{A} + \frac{p'}{p}\big) \psi_{z}'+ \frac{1}{p}\big( \frac{(pA')'}{A} + q+z \big) \psi_{z}=0.
\label{eq:Jo6}\end{equation}
By definitions \e{eq:Jo1}  and \e{eq:Jo3}  the coefficient at $\psi_{z}'$ here equals
\begin{equation}
 2\frac{A'}{A} + \frac{p'}{p} = 2\frac{a'}{a} + 2i\xi+ \frac{p'}{p}= 2i\xi  -\frac{\xi'}{\xi} .
\label{eq:Jo7}\end{equation}
Next, we compute the coefficient at $\psi_{z}$ in \e{eq:Jo6}.   Using \e{eq:Jo1A}, we see that
\[
 pA'= (pa' + i a^{-1}) e^{i\Xi}
\]
and
\[
 (pA')'= \big( (pa')'- a^{-1}\xi \big) e^{i\Xi}
\]
whence
\begin{equation}
A^{-1}(pA')'  + q
 =  a^{-1}(pa')' .
\label{eq:Jo8}\end{equation}
Substituting now expressions  \e{eq:Jo7}  and  \e{eq:Jo8} into  \e{eq:Jo6}, we obtain equation  \e{eq:Jo5}  with the coefficient
 \[
\rho_{z} (x) = (p  \xi )^{-1} \big(a^{-1} (pa')'  +z\big).
\] 
In view of  \e{eq:Jo1A} this coincides with definition \e{eq:Jr}.
    \end{pf}
    
    
     Next, we reduce the differential equation \e{eq:Jo5}  to a Volterra integral equation
       \begin{equation}
 \psi_{z} (x)=1+ (2i)^{-1}\int_{x}^\infty \big(1-  e^{-2i  \Xi(x)} e^{2i \Xi(y) } \big) \rho_{z} (y)\psi_{z} (y) dy ,\q x\geq x_{0}  .
  \label{eq:int}\end{equation} 
    
    \begin{lemma}\label{JOST2}
    Let assumptions  \e{eq:JJ}  be satisfied. Then equation \e{eq:int}  has a unique solution $ \psi_{z} (x)$. This solution satisfies  equation \e{eq:Jo5} and   
      \begin{equation}
 \psi_{z} (x)\to 1, \q \psi_{z}' (x)= o( \xi(x) ) \q \mbox{as}\q x\to\infty.
  \label{eq:int1a}\end{equation} 
       \end{lemma} 
    
    \begin{pf}
    Note that $\rho_{z}\in L^1 (x_{0},\infty)$ according to conditions \e{eq:JJ}. Therefore a bounded solution $\psi_{z}(x)$ of equation \e{eq:int}  can be standardly constructed by iterations.   
    
    Let us check that $\psi_{z}(x)$ satisfies equation \e{eq:Jo5}. Differentiating \e{eq:int}, we see that
    \begin{equation}
 \psi_{z}' (x)=  \xi(x) e^{- 2i \Xi(x) }  \int_{x}^\infty   e^{ 2i  \Xi(y)}   \rho_{z} (y)\psi_{z} (y) dy  
  \label{eq:int1}\end{equation} 
   and
     \begin{multline}
 \psi_{z}'' (x)= -  \xi(x)  \rho_{z} (x)\psi_{z} (x)  
 \\
 +  (\xi' (x) -2i  \xi^2(x) ) e^{-2i \Xi(x)}\int_{x}^\infty   e^{2i  \Xi(y)}   \rho_{z} (y)\psi_{z} (y) dy  .
  \label{eq:int2}\end{multline} 
  Substituting expressions  \e{eq:int1} and \e{eq:int2}  into the left-hand side of \e{eq:Jo5}, we see that it equals zero.
  
  Relations \e{eq:int1a} are direct consequences of \e{eq:int} and  \e{eq:int1}.
       \end{pf}
       
          Now it is easy    to conclude the proof of Theorem~\ref{JOST}. 
       
         \begin{pf}  
        Define the function  $ f_{z}  (x)$  by formula  \e{eq:Jo4}.  According to Lemma~\ref{JOST1} it satisfies    differential equation \e{eq:Jy}, and   according to Lemma~\ref{JOST2} it  has asymptotics \e{eq:Jo2}.  Moreover, differentiating \e{eq:Jo4}, we obtain  
  \[
f_{z}' (x)= i\xi(x) f_{z} (x) + a(x) e^{i\Xi (x)} \psi_{z}'(x) + a'(x) e^{i\Xi (x)} \psi_{z} (x).
\]
As we have already seen, the first term on the right has asymptotics \e{eq:Jo2+}. It follows from relations \e{eq:int1a} that 
the second and  third terms are $o(  \xi(x) a(x))$ and $o(a'(x))$, respectively. So, it remains to observe that
\[
o(a' (x))= o (p (x)^{-1} a (x)^{-1})=  o(  \xi(x) a(x))
\]
 according to condition \e{eq:JJ+} and identity \e{eq:Jo1A}.
    \end{pf}

 Theorem~\ref{JOST}  defines the Jost solutions for $x>x_{0}$.  Then the functions $f_{z} (x)$ are extended to all $x\geq 0$ as solutions of the 
differential equation \e{eq:Jy}.
    
    Let us introduce the second solution  $ \bar{f}_{\bar{z}}  (x)$ of differential equation \e{eq:Jy}. It follows from asymptotics  \e{eq:Jo2},  \e{eq:Jo2+} and equality \e{eq:Jo1A}  that the Wronskian
      \begin{equation}
\{ f_{z}  ,   \bar{f}_{\bar{z}}  \}  = 2 i \lim_{x\to \infty }\big( p(x) \xi (x) a^2 (x)\big) =2 i
  \label{eq:Jo1B}\end{equation}
 so that  these solutions  are linearly independent.
 
 We also observe that a solution $f_{z}(x)$ of equation \e{eq:Jy}  is determined uniquely by conditions \e{eq:Jo2} and  \e{eq:Jo2+}.  Indeed, if $\ti{f}_{z}(x)$ is another solution of equation \e{eq:Jy} satisfying these conditions, then  
 the Wronskian  $\{ f_{z}  ,   \ti{f}_z  \} =0$ so that  $\ti{f}_{z}(x)= c f_{z}(x)$ for some $c\in{\Bbb C}$. This constant equals $1$ according again to  \e{eq:Jo2}.
 
 Thus, the following result is a direct consequence of Theorem~\ref{JOST}.

 \begin{proposition}\label{JOLC}
 Under the  assumptions
  \e{eq:JJ} and \e{eq:JJ+} equation \e{eq:Jy} is in the LC case $($at infinity$)$.
     \end{proposition} 
     
     This result is not  really new; cf. Theorem~XIII.6.20 in the book \cite{DuSc}. 
    
      \subsection{Arbitrary solutions of the homogeneous equation}

An arbitrary solution   of the Schr\"odinger equation  \e{eq:Jy} is a linear combination of the Jost solutions $f_{z} (  x )$ and $ \bar{f}_{\bar{z}}(  x )$. In particular, this is true for regular solutions  $\varphi_{z} (  x )$ and 
$\theta_{z} (  x )$ distinguished by the boundary conditions  \e{eq:Pz} or \e{eq:Pzi}:
  \begin{equation}
\varphi_{z} (  x )= \sigma_{+} (z) f_{z}(  x ) + \sigma_{-} (z)    \bar{f}_{\bar{z}} (  x )
\label{eq:LC}\end{equation}
and
  \begin{equation}
\theta_{z} (  x )= \tau_{+} (z) f_{z} (  x ) + \tau_{-} (z)  \bar{f}_{\bar{z}} (  x ),
\label{eq:LCq}\end{equation}
where the coefficients $ \sigma_{\pm} (z)$ and $ \tau_{\pm} (z)$
 can be expressed via the Wronskians:
     \begin{equation}
       \begin{split}
  2i \sigma_{+}(z)=  \{\varphi_{z},   \bar{f}_{\bar{z}}  \},  \q
 2i \sigma_{-}(z)= -   \{\varphi_{z}, f _{z}  \}, 
\\
  2i \tau_{+}(z)=  \{\theta_{z},  \bar{f}_{\bar{z}}   \},  \q
 2i \tau_{-}(z)= -   \{\theta_{z}, f_{z}  \} .
   \end{split}
   \label{eq:LCx}\end{equation}

Observe that
 \begin{equation}
\sigma_{-} (z)= \ov{\sigma_{+}(\bar{z})}\q \mbox{and}\q \tau_{-} (z)= \ov{\tau_{+}(\bar{z})}
 \label{eq:LS}\end{equation}
 because $\varphi_{z}(x)= \ov{\varphi_{\bar{z}}(x)}$ and  $\theta_{z}(x)= \ov{\theta_{\bar{z}}(x)}$.
Of course, all coefficients $ \sigma_{\pm}(z)$ and $ \tau_{\pm}(z)$  are  entire functions of $z$.

According to \e{eq:LC}  and \e{eq:LCq}   the following result is a direct consequence of Theorem~\ref{JOST}. 

 \begin{theorem}\label{LC}
    Under the assumptions of Theorem~\ref{JOST} the solutions $\varphi_{z} (  x )$ and $\theta_{z} (  x )$ of equation \e{eq:Jy}
 have   asymptotics 
   \begin{equation}
\varphi_{z} (  x )=a (  x ) \big(\sigma_{+} (z)e^{ i \Xi (x)} + \sigma_{-}(z) e^{- i \Xi (x)}  + o(1)   \big)      
\label{eq:LC1P}\end{equation} 
and
 \begin{equation}
\theta_{z} (  x )= a (x) \big(\tau_{+} (z)e^{ i \Xi(x)} + \tau_{-} (z) e^{-i \Xi(x)}  + o(1)  \big)    
\label{eq:LC2q}\end{equation} 
 as $x\to\infty$.
These asymptotic formulas
can be differentiated in $x$;  in particular,
  \begin{equation}
\varphi_{z}' (  x )=i \xi(x) a (  x ) \big(\sigma_{+} (z)e^{ i \Xi (x)} - \sigma_{-}(z) e^{- i \Xi (x)}  + o(1)   \big)      .
\label{eq:LC1P+}\end{equation} 
   \end{theorem}

    In view of conditions \e{eq:Pz}  or \e{eq:Pzi}  the Wronskian $\{ \varphi_{z}, \theta_{z}\} =1$. On the other hand, we can calculate this Wronskian using relations   \e{eq:Jo1B} and \e{eq:LC}, \e{eq:LCq}. This yields an identity
   \begin{equation}
  2i  \big(\sigma_{+} (z) \tau_{-} (z) -\sigma_{-} (z)\tau_{+} (z)\big)=1,\q \forall z\in{\Bbb C}.
\label{eq:Wro}\end{equation}

Below, we need also the following fact.

\begin{proposition}\label{LCx}
    Under the assumptions of Theorem~\ref{JOST} we have an identity
     \begin{equation}
|\sigma_{-} (z)|^{2} - |\sigma_{+} (z)|^{2}  =\Im z    \int_{0}^{\infty}  |\varphi_{z} (x)|^{2} dx.
 \label{eq:LC2p}\end{equation} 
       \end{proposition}

\begin{pf}
Multiplying equation \e{eq:Jy} for $\varphi_{z}(x)$ by $\bar{\varphi}_{z}(x)$, integrating and taking the imaginary part, we see that
\[
- \Im \int_{0}^x \big( p(y) \varphi'_{z}(y) \big)'  \bar{\varphi}_{z}(y) dy=\Im z \int_{0}^x | \varphi_{z}(y) |^2 dy.
\]
Next, we integrate by parts on the left and take into  account the boundary condition \e{eq:BCz} whence
\begin{equation}
- \Im   \big( p(x) \varphi'_{z}(x)  \bar{\varphi}_{z}(x) \big)=\Im z \int_{0}^x | \varphi_{z}(y) |^2 dy.
\label{eq:WroW}\end{equation}
It follows from asymptotic formulas \e{eq:LC1P} and  \e{eq:LC1P+} that
  \begin{multline*}
 p(x) \varphi'_{z}(x)  \bar{\varphi}_{z}(x)  = i\big( |\sigma_{+} (z)|^{2} - |\sigma_{-} (z)|^{2}\big)
\\
 +i \big(\sigma_{+} (z) \bar{\sigma}_{-} (z) e^{2i \Xi (x)} -\sigma_{-} (z) \bar{\sigma}_{+} (z) e^{-2i \Xi (x)} \big) 
 + o(1) 
  \end{multline*}
where equality \e{eq:Jo1A} has been used.  Let us take the imaginary part of this expression. Then the second term on the right disappears. Substituting this expression into \e{eq:WroW} and passing to the limit $x\to\infty$, we arrive at \e{eq:LC2p}.
\end{pf}
 
  \subsection{Conditions on the coefficients}
  
  Let us discuss the assumptions of Theorem~\ref{JOST}.  The main condition is $ a\in L^2 ({\Bbb R}_{+})$. It requires that the product $p(x) q(x)\to\infty$ sufficiently rapidly, roughly speaking, faster than $x^2$.  The second 
  inclusion  \e{eq:JJ} as well as  
   condition \e{eq:JJ+} play  auxiliary roles. They mean that $p(x)$ does not grow too  rapidly compared to $a(x)$ and exclude too wild oscillations of the functions $p(x)$ and $a(x)$. For example, for the functions $p(x) = x^\beta$, $q(x)=x^\gamma$  (for large $x$) conditions \e{eq:JJ} and \e{eq:JJ+}  are satisfied if
 \begin{equation}
  \beta+\gamma>2 \q \mbox{and} \q \beta-\gamma <2.
  \label{eq:pa}\end{equation}
  This implies that, necessarily, $\gamma>0$, but it may be an arbitrary small number.  Observe  that $\beta\to 2$ if $\gamma\to 0$.
  
   It is noteworthy that  \e{eq:pa}  allows negative 
  $\beta$ provided  $  \gamma > 2+ |\beta|$.  In particular, according to Proposition~\ref{JOLC}  for such $\beta$ and $\gamma$ the operator  $H_{\rm min}$ is in the LC case. Condition on $\gamma$ is very important here. Indeed, the results of
    \cite{Y/HS} show  that if $q(x)=0$ and $p(x)\to 0$ very rapidly, then the corresponding Schr\"odinger operator is self-adjoint, its spectrum is absolutely continuous and coincides with $[0,\infty)$.


  \section{Schr\"odinger operators    and their quasiresolvents}
  
  We refer to the books \cite{Nai}, \S 17, and \cite{RS},  Sect.~X.1, for background information on the theory of symmetric differential operators.

  
    \subsection{Minimal and maximal operators}
    
    We here consider differential operators  \e{eq:ZP+} in  the space $L^2 ({\Bbb R}_{+})$. 
      The scalar product in  this space   is denoted $\la \cdot, \cdot\ra$;  $I$ is the identity operator. 
    
  Let us first define a minimal
    operator $H_{00}$    by the equality $H_{00} u= {\sf H} u$ on domain 
 $  {\cal D} (H_{00})  $ that consists of functions  $u\in C ^2 ({\Bbb R}_{+}) $ such that $u(x)=0$ for sufficiently large $x$,   limits $u(+0)=: u(0)$, $u'(+0)=: u'(0)$ exist  and condition \e{eq:BCz} is satisfied.
 Thus, the boundary condition \e{eq:BCz} at $x=0$ is included in the definition of the operator $H_{00}$ so that its self-adjoint  extensions are determined by conditions for $x\to\infty$ only.

 The closure   of  $H_{00}$ will be denoted  $H_{\rm min} $.  This operator is
symmetric  in the space $L^2 ({\Bbb R}_{+})$, and under  assumptions of this paper  its domain $  {\cal D} (H_{\rm min})  $ can be described efficiently (see Proposition~\ref{ABS}).  The adjoint operator  $H^*_{\rm min} =: H_{\rm max}$ is  again given  by the formula $H_{\rm max} u={\sf H} u$ on a set $  {\cal D} (H_{\rm max})$ consisting of functions $u$ in the local Sobolev space ${\Bbb H}^2_{\rm loc}$, satisfying  boundary condition \e{eq:BCz} and such that $u\in L^2 ({\Bbb R}_{+})$  and  ${\sf H} u\in L^2 ({\Bbb R}_{+})$.
In the LC case, the operator $H_{\rm max}$  is not   symmetric.  Integrating by parts, we see that
     for all $u, v \in {\cal D} (H_{\rm max})$ 
           \begin{equation}
\la {\sf H }u,v \ra - \la u,  {\sf H} v \ra = \lim_{x\to\infty}  p(x) \big(u(x)\bar{v}' (x)- u' (x) \bar{v}(x)\big)
\label{eq:qres4}\end{equation}
where the limit in the right-hand side exists but  is not necessarily zero.

 Recall that
\[
  H_{\min}= H_{\min}^{**}  = H_{\max}^{*}. 
  \]
  The operator  $H_{\rm min}$  is  self-adjoint if and only if the LP case occurs.
  In this paper we are interested in the LC case     when 
       \[
  H_{\rm min}\neq H_{\rm max}= H_{\rm min}^*.
 \]

Self-adjoint extensions $H$ of the operator $  H_{\rm min}$ satisfy the condition
  \[
H_{\rm min}\subset H=H^*\subset H_{\rm min}^*=: H_{\rm max}.
  \]
In the LC case  the operators $H_{\rm max}$ are not symmetric.

Since the operator $H_{\rm min}  $  commutes with the complex conjugation,  its deficiency indices 
   \[
   d_{\pm}: =\dim\ker (H_{\rm max}-z I),\q \pm \Im z>0,
   \]
     are equal, i.e. $d_{+} =d_{-}=:d $, and, so, $H_{\rm min}  $ admits self-adjoint extensions.   For an arbitrary $z\in{\Bbb C}$, all solutions  of equation \e{eq:Jy}  with boundary condition 
     \e{eq:BCz}
       are given by the formula $u (x)= \Gamma \varphi_{z}(x)$  for some $\Gamma \in{\Bbb C}$. They belong to $ {\cal D} (H_{\rm max})$ if and only if $\varphi_{z} \in  L^2 ({\Bbb R}_{+})$.  Therefore $d=0$ if  $\varphi_{z} \not\in  L^2 ({\Bbb R}_{+})$ for $\Im z\neq 0$; otherwise $d=1$.

      \subsection{Quasiresolvent of the maximal operator}

Recall that in the LC  case
 inclusions \e{eq:PQz}  are satisfied.    Following \cite{Schr-LC}, let us define, for all $z\in {\Bbb C}$,  a bounded   operator ${\cal R}  (z)$ in the  space $L^2 ({\Bbb R}_{+})$  by   the equality 
   \begin{equation}
 ( {\cal R} (z)h) (x) = \theta_{z} (x) \int_{0}^x \varphi_{z} (y) h(y) dy+ 
\varphi_{z} (x) \int_x^\infty \theta_{z} (y) h(y) dy . 
 \label{eq:RR11}\end{equation} 
  Actually, the operator ${\cal R} (z)$ belongs to the Hilbert-Schmidt class. It depends analytically on $z\in {\Bbb C}$ and ${\cal R} (z)^*={\cal R} (\bar{z})$. 
We prove (see  Theorem~\ref{res}) that, in  a natural sense, 
${\cal R} (z)$ plays the role of the resolvent of the operator $H_{\rm max}$.  We call it
the  quasiresolvent of the operator $H_{\rm max}$.   It

 Let us enumerate some simple properties of the operator ${\cal R}  (z)$.  Differentiating  
   definition  \e{eq:RR11}, we see that 
    \begin{equation}
 ( {\cal R} (z)h)' (x) =   \theta_{z}' (x) \int_{0}^x \varphi_{z} (y) h(y) dy+ 
 \varphi_{z}' (x) \int_x^\infty \theta_{z} (y) h(y) dy 
 \label{eq:RR12}\end{equation} 
 for all $h\in L^2 ({\Bbb R}_{+})$.
 In particular, it follows from relations  \e{eq:RR11} and \e{eq:RR12}  that
  \begin{equation}
({\cal R}  (z)h)(0)= \varphi_{z} (0)\la h, \theta_{\bar{z}}\ra 
\label{eq:r01}\end{equation}
and
 \begin{equation}
({\cal R}  (z)h)'(0)= \varphi_{z}' (0)\la h, \theta_{\bar{z}}\ra 
\label{eq:r02}\end{equation}
where $\varphi_{z} (0)$ and $\varphi_{z}' (0)$  are defined by equalities \e{eq:Pz}  or \e{eq:Pzi}.

A proof of the following statement is  close to the construction  of the resolvent for  essentially self-adjoint Schr\"odinger operators.    
   
    \begin{theorem}\label{res}
    Let  inclusions \e{eq:PQz}  hold  true.
 For all $z \in {\Bbb C}$, we have
  \begin{equation}
{\cal R}  (z): L^2 ({\Bbb R}_{+})\to  {\cal D} (H_{\rm max} )
\label{eq:qres}\end{equation}
and
  \begin{equation}
(H_{\rm max} -zI) {\cal R}  (z)=I.
\label{eq:qres1}\end{equation}
   \end{theorem} 
 
  \begin{pf} 
  Let $h\in L^2 ({\Bbb R}_{+})$ and $u(x)=    ({\cal R} (z)h)(x)$.   Boundary condition \e{eq:BCz}  is a direct consequence of relations  \e{eq:r01}  and  \e{eq:r02}. Differentiating \e{eq:RR12}, we see that
  \begin{multline}
 (p(x)u' (x))' =   (p(x) \theta_{z}' (x) )'\int_{0}^x \varphi_{z} (y) h(y) dy
 \\
 + (p(x) \varphi_{z}' (x) )' \int_x^\infty \theta_{z} (y) h(y) dy + p(x) \big( \theta'_{z}(x) \varphi_{z}(x)-\theta
 _{z}(x) \varphi'_{z}(x)\big)  h(x).
 \label{eq:RR13}\end{multline} 
Since the Wronskian $\{ \varphi_{z}, \theta_{z}\}=1 $,   the last term in the right-hand side equals $-h(x)$.  Putting now together  together equalities \e{eq:RR11}  and \e{eq:RR13} and using equation \e{eq:Jy} for the functions $\varphi_{z}(x)$ and $\theta_{z}(x)$, we obtain   the equation 
\[
- (p(x) u' (x))'  + q(x) u(x)-z u(x)=h(x) 
  \]
  where   $h\in L^2 ({\Bbb R}_{+})$. 
  Together with boundary condition \e{eq:BCz} this implies that $   H_{\rm max} u -z u=h $.  This yields both \e{eq:qres} and \e{eq:qres1}. 
    \end{pf}


    Note that  solutions $u(x)$ of differential equation   \e{eq:Jy}  satisfying condition \e{eq:BCz} are given by the formula $u(x)= \Gamma \varphi_{z} (x)$  for some $ \Gamma \in{\Bbb C}$. Therefore we
    can   state 
    
     \begin{corollary}\label{res1}
     All solutions of the equation
       \[
(H_{\rm max} -zI) u =h \q \mbox{where}  \q z\in {\Bbb C} \q \mbox{and}\q h\in L^2 ({\Bbb R}_{+})
\]
  for $u\in {\cal D} (H_{\rm max} )$ are  given by the formula
 \begin{equation}
  u = \Gamma \varphi_{z}+ {\cal R}  (z) h \q \mbox{for some} \q \Gamma=\Gamma (z;h) \in{\Bbb C}.
\label{eq:qres3}\end{equation}
   \end{corollary}

An  asymptotic relation for $ ( {\cal R} (z)h) (x)$ is a direct consequence of definition  \e{eq:RR11}
   and condition  \e{eq:PQz}:
      \begin{equation}
 ( {\cal R} (z)h) (x) =  \theta_{z} (x) \la    h, \varphi_{\bar{z}}  \ra+   o(|\varphi_{z} (x)| +|\theta_{z} (x)| )\q {\rm as}\q x\to\infty.
 \label{eq:Ras}\end{equation}

 \section{  Self-adjoint extensions and their resolvents }
 
 Here we find an asymptotic behavior as $|x|\to\infty$ of all functions $u(x)$ in the domain of the maximal operator
 $H_{\rm max}$. This allows us to give an efficient description of all self-adjoint extensions of the operator
  $H_{\rm min}$.

  \subsection{Domains of maximal operators}
  
  Recall that  boundary condition \e{eq:BCz} at $x=0$ is included in our  definition of the minimal operator
 $H_{\rm min}$. Our goal now is to find a similar condition for $x\to\infty$ distinguishing self-adjoint extensions of 
  $H_{\rm min}$.  
  
   The starting point of our construction is  asymptotic formula \e{eq:LC2}  for   functions  $u\in H_{\rm max}$. 
    Recall that the amplitude $a(x)$ and the phase $\Xi(x)$ were defined by formulas \e{eq:Jo1}. The coefficients $\sigma_{\pm} (z)$, $\tau_{\pm}(z)$ are given by equalities \e{eq:LCx}, and the number $\Gamma (z; h)$ is determined by formula \e{eq:qres3}.

  \begin{theorem}\label{LCR}
  Let the assumptions of Theorem~\ref{JOST} be satisfied.  Then an arbitrary function
   $u\in{\cal D}(H_{\rm max})$ 
 has asymptotics \e{eq:LC2}
with some   coefficients $s_{\pm} = s_{\pm} (u) $.     They can be constructed by relations
 \begin{equation}
  \begin{split}
s_{+}(u)=  \Gamma(z; ({\sf H}-zI)u) \sigma_{+} (z)+ \la ( {\sf H}-zI) u, \varphi_{\bar{z}}\ra \tau_{+} (z),
  \\
 s_{-}(u)= \Gamma(z; ({\sf H} -zI)u) \sigma_{-} (z)+ \la ({\sf H} -zI) u, \varphi_{\bar{z}}\ra \tau_{-} (z)
  \end{split}
  \label{eq:LC31}\end{equation}
  where the number   $z\in {\Bbb C}$ is arbitrary.

  Conversely, for arbitrary $ s_{+},s_{-}  \in {\Bbb C} $, there exists a  function 
$u \in{\cal D}(J_{\rm max})$ such that    asymptotics \e{eq:LC2}  holds.
 \end{theorem}
 
 
   \begin{pf}
   According to Corollary~\ref{res1}  a function    $u \in{\cal D}(H_{\rm max})$  admits representation  \e{eq:qres3} where the operator ${\cal R} (z)$  is defined by equality \e{eq:RR11}.
      In view of   relation \e{eq:Ras}  and asymptotics \e{eq:LC1P}, \e{eq:LC2q} we have
    \begin{equation}
 ( {\cal R} (z)h) (x) =    a (x)  \big( \tau_{+} (z)e^{i \Xi (x) }+ \tau_{-} (z) e^{-i \Xi (x)} \big) \la h, \varphi_{\bar{z} }\ra     + o(a (x) ), \q x\to\infty,
\label{eq:RR1r}\end{equation}
  for all functions $h \in L^2 ({\Bbb R}_{+})$.  Therefore it follows from \e{eq:qres3}  that
    \begin{multline*}
u (x) = a (x) \Gamma(z; ({\sf H}-zI)u)  \big(\sigma_{+} (z) e^{i \Xi(x)} + \sigma_{-} (z)e^{-i \Xi(x)}\big) 
\\
+ a (x) \big(\tau_{+} (z) e^{i \Xi(x)} + \tau_{-} (z)e^{-i \Xi(x)}\big) \la  ( {\sf H}-zI) u, \varphi_{\bar{z} }\ra    + o(a (x))   
 \end{multline*} 
as $ x\to\infty$. This yields  relation \e{eq:LC2} with the coefficients $s_{\pm}$ defined by \e{eq:LC31}.

Conversely, given $ s_{+} $ and $ s_{-} $ and fixing some $z\in{\Bbb C}$, we consider   a system of equations 
 \begin{equation}
  \begin{split}
s_{+}=  \Gamma\sigma_{+} (z)+ \la h , \varphi_{\bar{z} }\ra \tau_{+} (z),
  \\
 s_{-} = \Gamma \sigma_{-} (z)+ \la h , \varphi_{\bar{z} }\ra \tau_{-} (z).
  \end{split}
  \label{eq:LC3x}\end{equation}  
   for  $\Gamma$ and $ \la h, \varphi_{\bar{z} }\ra $. According to \e{eq:Wro} the determinant  of this system is not zero so that  $\Gamma$ and $ \la h, \varphi_{\bar{z} }\ra $ are uniquely determined by  $s_{+} $ and $ s_{-} $.
 Then we take any $h$ such that its scalar product with $ \varphi_{\bar{z} }$ equals the found value of  $ \la h, \varphi_{\bar{z} }\ra$.
  Finally, we define $u$ by formula  \e{eq:qres3}. Asymptotics as $x\to\infty$ of $\varphi_{z} (x)$   and $( {\cal R} (z)h) (x)$ are given by formulas \e{eq:LC1P} and \e{eq:RR1r}, respectively. In view of equations \e{eq:LC3x} this leads to asymptotics \e{eq:LC2}.
       \end{pf}

       
      Theorem~\ref{LCR}  yields a mapping $ {\cal D}(H_{\rm max})\to{\Bbb C}^2$ defined by the formula
       \begin{equation}
u\mapsto (s_{+}  (u) , s_{-}(u)).
\label{eq:mapping}\end{equation}
      The construction of      Theorem~\ref{LCR} depends on the choice of $z\in{\Bbb C}$, but this mapping is defined intrinsically.  In particular, we can set $z=0$ in all formulas of  Theorem~\ref{LCR}.  Note that mapping \e{eq:mapping} is surjective.  We also observe that \e{eq:mapping} 
        plays the role of a mapping $u\mapsto(u(0), u' (0))$   for functions in  the local Sobolev class ${\Bbb H}_{\rm loc}^2$.
        

     Under the assumptions of Theorem~\ref{JOST} the right-hand side of \e{eq:qres4} can be expressed  in terms of the coefficients $s_{+}$ and $s_{-}$.

 \begin{proposition}\label{resAS}
     For all $u, v \in {\cal D} (H_{\rm max})$, we have an identity
    \begin{equation}
\la H_{\rm max}u,v \ra - \la u, H_{\rm max}v \ra = 2 i \big( s_{+}(u) \ov{s_{+}(v) }- s_{-}(u) \ov{s_{-}(v) }\big)  .
\label{eq:ABS}\end{equation}
 \end{proposition} 
   
    \begin{pf}
    Let us proceed from equality  \e{eq:qres4}.
    It follows from formula \e{eq:LC2}  that
       \begin{multline*}
- i p(x) u'(x) \bar{v} (x)  =    \big(s_{+}(u)   e^{ i \Xi (x)} - s_{-} (u)  e^{- i \Xi (x)}\big)\big(\ov{s_{+}(v) }  e^{- i \Xi (x)} +  \ov{ s_{-} (v)} e^{ i \Xi (x)}\big) + o(1)
\\
=  s_{+}(u) \ov{s_{+}(v) }-  s_{-}(u) \ov{s_{-}(v) } + s_{+}(u) \ov{s_{-}(v) } e^{2 i \Xi (x)} 
-s_{-}(u) \ov{s_{+}(v) } e^{-2 i \Xi (x)}+ o(1).
 \end{multline*}
 and, similarly,
 \begin{multline*}
 i p(x) u(x) \bar{v}' (x)  =    \big(s_{+}(u)   e^{ i \Xi (x)} + s_{-} (u)  e^{- i \Xi (x)}\big)
 \big(\ov{s_{+}(v) }  e^{- i \Xi (x)} -  \ov{ s_{-} (v)} e^{ i \Xi (x)}\big) + o(1)
\\
=  s_{+}(u) \ov{s_{+}(v) }-  s_{-}(u) \ov{s_{-}(v) } - s_{+}(u) \ov{s_{-}(v) } e^{2 i \Xi (x)} 
+s_{-}(u) \ov{s_{+}(v) } e^{-2 i \Xi (x)}+ o(1).
 \end{multline*}
 Let us take the sum of the last two expressions and observe that the terms containing $e^{2 i \Xi (x)}$ and $e^{-2 i \Xi (x)}$ cancel each other.  This yields
    \[
- i p(x) (u'(x) \bar{v} (x) - u (x)\bar{v}' (x))=  2 s_{+}(u) \ov{s_{+}(v) }- 2 s_{-}(u) \ov{s_{-}(v) } + o(1).
 \]
Passing here to the limit $n\to\infty$ and using equality \e{eq:qres4}, we obtain identity \e{eq:ABS}.
       \end{pf}

    We can now characterize the set     $ {\cal D} (H_{\rm min})$.
       
       
       \begin{proposition}\label{ABS}
       A vector
    $v  \in {\cal D} (H_{\rm max})$ belongs to     $ {\cal D} (H_{\rm min})$ if and only if
    $v (x)= o(a (x))$ as $x\to\infty$, that is, 
     \begin{equation}
    s_{+} (v) = s_{-} (v) =0.
      \label{eq:ABS4}\end{equation}
    \end{proposition} 
     
     \begin{pf}    
     A vector $v$ belongs to     $ {\cal D} ( H_{\max}^*)$ if and only if
       \begin{equation}
    \la H_{\max}u,v \ra= \la u, H_{\max}v \ra
  \label{eq:ABS2}\end{equation}
    for all   $u  \in {\cal D} (H_{\rm max})$.  According to Proposition~\ref{resAS}  equality \e{eq:ABS2} is equivalent to
     \begin{equation}
 s_{+}(u) \ov{s_{+}(v) }- s_{-}(u) \ov{s_{-}(v) }=0.
  \label{eq:ABS3}\end{equation}
  This  is of course true if \e{eq:ABS4} is satisfied. Conversely, if  \e{eq:ABS3} is satisfied for all  $u\in {\cal D} ( H_{\max})$, we use that according to Theorem~\ref{LCR}  the numbers $ s_{+}(u)$ and $s_{-}(u)$  are arbitrary. This implies \e{eq:ABS4}.
         \end{pf}       
    
This result shows that  \e{eq:mapping} considered as a mapping of the factor space  $  {\cal D} (H_{\rm max})/ 
  {\cal D} (H_{\rm min})$ onto ${\Bbb C}^2$ is injective.
  
    \subsection{Self-adjoint extensions}
    
    All self-adjoint extensions $ H_{\omega} $ of the operator $H_{\min}$ are parametrized by complex numbers $\omega\in{\Bbb T}\subset {\Bbb C}$.  Let a set $ {\cal D} (  H_\omega)\subset {\cal D} (    H_{\max})$ of vectors $u$  be distinguished by  condition     \e{eq:gen2}.

      \begin{theorem}\label{ABX}
           Let the assumptions of Theorem~\ref{JOST}  be satisfied. Then
     for all $\omega\in{\Bbb T}$, the operators $H_{\omega}$ are  self-adjoint. Conversely, every   operator $H$ such that 
     \begin{equation}
  H_{\min} \subset  H =H^*\subset H_{\max} 
      \label{eq:ABX2}\end{equation}
      equals $H_{\omega}$ for some $\omega\in{\Bbb T}$.
    \end{theorem} 
    
    \begin{pf}
    We proceed from Proposition~\ref{resAS}.   If $u,v \in  {\cal D} (  H_\omega)$,  it follows from condition \e{eq:gen2} that   $s_{+}(u) \ov{s_{+}(v) }= s_{-}(u) \ov{s_{-}(v) }$. Therefore   according to equality \e{eq:ABS} $\la H_{\omega}u,v  \ra= \la u, H_{\omega}v \ra $ whence $H_{\omega}\subset H_{\omega}^*$. If $ v\in {\cal D} (H_\omega^*)$, then 
   $\la H_{\omega}u,v \ra= \la u, H_{\omega}^* v \ra $
    for all $ u\in {\cal D} ( H_\omega)$ so that in view of \e{eq:ABS} equality \e{eq:ABS3} is satisfied.
    Therefore $s_{-}(u) (\omega\ov{s_{+}(v) }- \ov{s_{-}(v) })=0$. Since $s_{-}(u)$ is arbitrary, we see  that  $ \omega\ov{s_{+}(v) }- \ov{s_{-}(v) }=0$, and hence $ v\in {\cal D} (H_\omega)$.
    
    Suppose that an operator $H$  satisfies conditions  \e{eq:ABX2}.  Since  $H$ is symmetric, it follows from  Proposition~\ref{resAS} that  equality   \e{eq:ABS3} equality is true  for all $ u, v\in {\cal D} (  H )$.  Setting here $u=v$, we see that $|  s_{+}(v)|=|   s_{-}(v)|$.   There exists a vector $v_{0} \in {\cal D} (  H )$ such that
    $s_{-}(v_{0})\neq 0$ because  $H\neq H_{\rm min}$. Let us set $\omega= s_{+}(v_{0})/ s_{-}(v_{0})$. Then $|\omega|=1$  and relation  \e{eq:gen2}  is a direct consequence of \e{eq:ABS3}.
        \end{pf}

       \subsection{Resolvent}
     
     Now it easy to construct the   resolvent  of the operator $H_\omega$  defined in the previous subsection.
     We previously note that, by definition \e{eq:LC2}, 
     \[
     s_{\pm} (\varphi_{z})= \sigma_{\pm} (z),  \q      s_{\pm} (\theta_{z})= \tau_{\pm} (z)
     \]
     and $|\sigma_{+} (z)|\neq |\sigma_{-} (z)|$ for $\Im z\neq 0$, by Proposition~\ref{LCx}.
     
     \begin{theorem}\label{ABX1}
            Let the assumptions of Theorem~\ref{JOST}   be satisfied. Then
     for all $z\in {\Bbb C}$ with $\Im z\neq 0$ and all $h\in L^2 ({\Bbb R}_{+})$, the resolvent $R_{\omega} (z)= (H_{\omega}-zI)^{-1}$ of the operator $H_\omega$ is given by   the equality 
      \begin{equation}
R_{\omega} (z) h =  {\gamma}_{\omega} (z) \la h,  \varphi_{\bar{z}}\ra \varphi_{z}+ {\cal R}(z)h
      \label{eq:ABY1}\end{equation}
      where
           \begin{equation}
 \gamma_{\omega}(z)= -\frac{\tau_{+}(z)-\omega \tau_{-}(z)} {\sigma_{+}(z)-\omega \sigma_{-}(z)}.
    \label{eq:gam}\end{equation}
    \end{theorem} 
    
         \begin{pf}
           According to      Corollary~\ref{res1} a vector  $u= R_{\omega}(z) h$ is given by formula \e{eq:qres3} where the coefficient $\Gamma$ is determined by condition 
            \e{eq:gen2}. It follows from Theorem~\ref{LCR}  than the function $u (x)$  has  asymptotics \e{eq:LC2}  with the coefficients $s_{\pm}$ defined by relations 
            \e{eq:LC3x}.  Thus, 
          $u\in {\cal D}( H_{\omega})$ if and only if  
    \[
    \Gamma \sigma_{+} (z) + \tau_{+} (z) \la h,  \varphi_{\bar{z}}\ra=\omega \big( \Gamma \sigma_{-} (z) + \tau_{-} (z)\la h,  \varphi_{\bar{z}}\ra\big)
    \]
    whence
    \[
\Gamma= -\frac{\tau_{+}(z)-\omega \tau_{-}(z)} {\sigma_{+}(z)-\omega \sigma_{-}(z)}\: \la h,  \varphi_{\bar{z}}\ra.
      \]
Substituting this expression into \e{eq:qres3}, we arrive at  formulas
    \e{eq:ABY1},   \e{eq:gam}.
          \end{pf}
          
          \begin{corollary}\label{abx}
          The resolvents $R_{\omega} (z)$ belong to the Hilbert-Schmidt class if $\Im z\neq 0$, whence the spectra of the operators $H_{\omega}$ are discrete.
    \end{corollary} 
    
    This result is   well known; see, e.g., Theorem~10 in Chapter~VII of the book \cite{Nai} or Theorem~5.8 in   the book \cite{Titch}.

          It follows from formula \e{eq:gam} that the spectra of the operator $H_{\omega}$ consist of the points  $z$ where
          \begin{equation}
          \sigma_{+}(z)-\omega \sigma_{-}(z)=0.
         \label{eq:ABY3}\end{equation}
                 Since the functions $\sigma_{+} (z)$ and $\sigma_{-} (z)$  are analytic, this again implies that 
                  the spectra of   $H_{\omega}$ are discrete. Of course the roots $z$ of equation \e{eq:ABY3} 
          lie        on the real axis because $|\sigma_{+} (z)|\neq |\sigma_{-} (z)|$ for $\Im z\neq 0$.   This fact  had  of course to be expected since $z$ are eigenvalues of the self-adjoint operator $H_{\omega}$.  We finally note that the discreteness of the spectrum of the operators $H_{\omega}$ is quite natural because their domains ${\cal D}( H_{\omega})$ are distinguished by boundary conditions at both ends of   ${\Bbb R}_{+}$. 
 Therefore  $H_{\omega}$ acquire some features of regular operators.

   \subsection{Spectral measure}
   
   In view of the spectral theorem, Theorem~\ref{ABX1} yields a representation  for   the Cauchy-Stieltjes transform of the spectral measure  $dE_\omega (\lambda)$ of the operator $H_{\omega}$.
   
   \begin{theorem}\label{RESc}
    Let  inclusions \e{eq:PQz}  hold.
     Then for all $z\in {\Bbb C}$ with $\Im z\neq 0$ and all $h\in L^2 ({\Bbb R}_{+})$,  we have  an equality
      \begin{equation}
     \int_{-\infty}^\infty
 (\lambda-z)^{-1} d (E_\omega (\lambda)h,h)= \gamma_{\omega} (z)|\la \varphi_{z},h\ra |^2 + ({\cal R}(z) h,h).
     \label{eq:E1}\end{equation}
    \end{theorem} 
    
    Recall that the operators ${\cal R}(z)$ are defined by formula \e{eq:RR11}. Therefore $({\cal R}(z) h,h)$ are entire functions of $z\in{\Bbb C}$,  and the singularities of the integral in \e{eq:E1} are determined by the function $\gamma_\omega (z)$.  Thus,   \e{eq:E1} can be considered as a modification of the classical    Nevanlinna formula (see his original paper  \cite{Nevan} or, for example, formula (7.6) in the book \cite{Schm}) for   the Cauchy-Stieltjes transform of the spectral measure in the  theory of Jacobi operators.  We mention however that, for Jacobi operators acting in the space $L^2 ({\Bbb Z}_{+})$, there is the canonical choice of a generating vector and  of a spectral measure. This is not the case for differential operators in $L^2 ({\Bbb R}_{+})$.

     We finally note an obvious fact: if $\lambda   $ is an eigenvalue of an operator $H_{\omega}$, then  corresponding eigenfunctions equal $c\varphi_{\lambda} (x)$ where $c\in{\Bbb C}$.

         \subsection{Concluding remarks}
         
         In the LC case, self-adjoint extensions of  the operator $H_{\rm min}$ are traditionally described by the following procedure; see the classical books \cite{CoLe}, Chapter~IX,  \cite{Nai}, \S 17, 18,  and the recent monograph  \cite{Schm}, Sect.~14.4, 15.3 and 15.4.  Take any real      functions $\varrho_j \in {\cal D} (H_{\rm max})$, $j=1,2$, such that
  \[
         \lim_{x\to\infty}  p(x) \big(\varrho_{1}' (x)\varrho_{2} (x)-\varrho_{1} (x)\varrho_{2}' (x)\big)=0,
         \]
         and set $\varkappa_{s} (x)= s \varrho_{1} (x)+ \varrho_2 (x)$ where $s\in{\Bbb R}$, $\varkappa _{\infty} (x)=  \varrho_{1} (x) $.  Let a set ${\cal D}_{s}\subset{\cal D} (H_{\rm max})$ be distinguished by the condition 
         \[
                \lim_{x\to\infty}  p(x) \big(u ' (x)\varkappa_{s} (x)- u  (x)\varkappa_{s}' (x)\big)=0 \q \mbox{for}\q u\in {\cal D}_{s},
                \]
                and let $\wt{H}_{s}$ be the restriction of $H_{\rm max}$ on the set 
                ${\cal D}_{s}=: {\cal D} (\wt{H}_{s})$.  Then the operators $\wt{H}_{s}$ are self-adjoint, and all self-adjoint extensions of the operator $H_{\rm min}$  coincide with one of the operators $\wt{H}_{s}$  for some $s\in {\Bbb R}\cup\{\infty\}$.  This construction does not look very efficient, in particular, because it depends on the choice of the functions $\varrho_1, \varrho_2$.
                
                Another possibility is to use von Neumann formulas. They were conveniently adapted to operators commuting with the complex conjugation in the survey \cite{Simon}, Theorem~2.6, and then applied in this paper to Jacobi operators (see also Sect.~16.3 in \cite{Schm}). Following \cite{Schr-LC}, we briefly describe  here this construction for Schr\"odinger operators \e{eq:ZP+}.  Recall that $\varphi_{0} (x)$ and $\theta_{0} (x)$ are the solutions of equation \e{eq:Jy}  where $z=0$ satisfying conditions \e{eq:Pz} or \e{eq:Pzi}.  We set $\ti{\theta}_{0} (x)=\omega (x) \theta_{0} (x)$ where  $\omega (x)$ is a smooth  function such that $\omega (x)=0$  for small $x$ and $\omega (x)=1$  for large $x$. Define operators $H^{(t)}$ as the restrictions of $H_{\rm max}$  on direct sums
            \begin{equation}
              {\cal D}(H_{\rm min} )\dotplus \{t\varphi_{0}+   \ti{\theta}_{0}\}=: {\cal D} (H^{(t)})\q \mbox{for} \q t\in {\Bbb R}
                \label{eq:EX}\end{equation}  
              and ${\cal D}(H_{\rm min}) \dotplus\{\varphi_{0}\}=: {\cal D} (H^{(\infty)})$.
        Then the operators $H^{(t)}$ are self-adjoint, and all self-adjoint extensions of the operator $H_{\rm min}$  coincide with one of the operators $H^{(t)}$ for some $t\in {\Bbb R}\cup\{\infty\}$. A drawback of this construction is that, apparently, the set $   {\cal D}(H_{\rm min} )$ cannot be described efficiently without some assumptions on the coefficients $p(x)$ and $q(x)$.
        
         Finally, we indicate a link between the operators $H^{(t)}$ and $H_{\omega}$ considered in this paper.   Let $u\in {\cal D} (H^{(t)})$
         where  $t\in {\Bbb R}$. According to Theorem~\ref{LC} and Proposition~\ref{ABS}  it follows from definition \e{eq:EX} that $u(x)$ has asymptotics \e{eq:LC2}  where 
         \[
         s_{+}(u)= t \sigma_{+}(0)+ \tau_{+}(0) \q\mbox{and}\q
           s_{-}(u)= t \sigma_{-}(0)+ \tau_{-}(0).
           \]
           Therefore relation  \e{eq:gen2}  is satisfied with
           \begin{equation}
           \omega=\frac{ t \sigma_{+}(0)+ \tau_{+}(0)}{ t \sigma_{-}(0)+ \tau_{-}(0)}.
                 \label{eq:EX1}\end{equation} 
                 If $u\in {\cal D} (H^{(\infty)})$, then  this equality holds true with $t=\infty$, that is, $\omega= \sigma_{+}(0)\sigma_{-}(0)^{-1}$.  Note that $|\omega|=1$ by virtue of identity \e{eq:LS}.  This shows that  $u\in {\cal D} (H_{\omega})$; see the definition at the beginning of Sect.~4.2. Conversely, if $u\in {\cal D} (H_{\omega})$, then $u\in {\cal D} (H^{(t)})$ with $t$  determined by \e{eq:EX1}.


\begin{thebibliography}{99}
    
 


 
   
   
    
  
  
  

    
 \bibitem {CoLe} E. A. Coddington and N. Levinson, {\it  Theory of ordinary differential equations},
McGraw-Hill, New York, 1955.

 \bibitem {DuSc} N. Dunford and J. T. Schwartz, 
 {\it Linear operators},   part  2, Interscience Publishers,
New York, London, Sydney, 1963.
 
       
       
 

 
   
  
  

 

 
    \bibitem{Nai} M.~A.~Naimark, \emph{Linear differential operators}, Ungar, New York,   1968.

  
   
  
  \bibitem {Nevan} R.~Nevanlinna, Asymptotische Entwickelungen beschr\"ankter Funktionen und das Stieltjessche  Momentenproblem, Ann. Acad. Sci. Fenn. A {\bf 18}, No. 5  (1922), 52 pp.
  
   
  
     \bibitem {Olver}  F.  W.  J.~Olver, {\it Introduction to asymptotics and special functions}, Academic Press, 1974.
 
       

  \bibitem {RS} M. Reed and B. Simon, {\it Methods of Modern Mathematical Physics} II, Academic Press, 1975.
  

 \bibitem{Schm} K.~Schm\"udgen, \emph{Unbounded self-adjoint operators on Hilbert space}, Springer, Dordrecht, Heidelberg, New York, London, 2012.
 
  \bibitem{Simon}
  B.~Simon, 
 The classical moment problem as a self-adjoint finite difference operator, Advances in  Math.  
{\bf 137} (1998), 82-203.
   

  
  

  

    
  
\bibitem {Titch}  E. C.~Titchmarsh, {\em Eigenfunction expansions associated with second-order differential equations}, vol. 1, Oxford, 1946.



  \bibitem{Y/HS} D. R. Yafaev, Spectral and scattering theory for differential and Hankel operators,     Advances in Math.   {\bf 308} (2017), 713-766.


 



 
  
    

          
         \bibitem{Jacobi-LC} D. R. Yafaev, Self-adjoint Jacobi operators in the limit circle case, arXiv  2104.13609.
         
          \bibitem{Schr-LC} D. R. Yafaev, Self-adjoint  differential operators in the limit circle case, arXiv  2105.08641.
         
         	 
  
            \end{thebibliography}
 \end{document}